\documentclass[11pt]{article}
\usepackage[leqno]{amsmath}
\usepackage{amsfonts,amssymb,theorem,amscd}
\usepackage[mathscr]{eucal}
\usepackage[leqno]{amsmath}
\usepackage{amsfonts,amssymb,theorem,amscd}
\usepackage[mathscr]{eucal}
      \makeatletter
\def\setseccntfmt{\renewcommand{\@seccntformat}[1]{\S
    \csname the##1\endcsname.\hspace{1ex}}}
\def\setsubseccntfmt{\renewcommand{\@seccntformat}[1]{%
    (\csname the##1\endcsname)\hspace{0.5ex}}}
\renewcommand{\section}{\setseccntfmt\@startsection
  {section}{1}{0mm}{-\baselineskip}{0.5\baselineskip}{\sf\bfseries\Large}}
\def\presubsection{\setsubseccntfmt\@startsection
  {subsection}{2}{0mm}{-\baselineskip}{-0.5ex}{\bfseries\upshape}}
\def\tmpa{}\def\tmpb{}
\newcommand{\addspaceifnonempty}[1]{\def\tmpa{}\def\tmpb{#1}%
  \ifx\tmpa\tmpb{}\else{\hspace{0.5ex}#1\hspace{1ex}}\fi}
\renewcommand{\subsection}[1][]{\presubsection{\addspaceifnonempty{#1}}}
\def\presubsubsection{\setsubseccntfmt\@startsection
  {subsubsection}{3}{0mm}{-\baselineskip}{-0.5ex}{\bfseries\upshape}}
\renewcommand{\subsubsection}[1][]{\presubsubsection{\addspaceifnonempty{#1}}}
\def\preparagraph{\setsubseccntfmt\@startsection
  {paragraph}{4}{0mm}{-\baselineskip}{-0.5ex}{\bfseries\upshape}}
\renewcommand{\paragraph}[1][]{\preparagraph{\addspaceifnonempty{#1}}}
\gdef\th@nonumplain{\normalfont\itshape
  \def\@begintheorem##1##2{\item[\hskip\labelsep\theorem@headerfont ##1]}%
  \def\@opargbegintheorem##1##2##3{%
    \item[\hskip\labelsep \theorem@headerfont ##1\ (##3)]}}
\gdef\th@change{%\normalfont\slshape
  \def\@begintheorem##1##2{\item[\hskip\labelsep
    {\bfseries\upshape(##2)\hskip 1ex}\theorem@headerfont ##1]}%
  \def\@opargbegintheorem##1##2##3{\item[\hskip\labelsep
    {\bfseries\upshape(##2)\hskip 1ex}\theorem@headerfont ##1\ (##3)]}}
\makeatother   \def\nthm{\newtheorem}
\theoremheaderfont{\normalfont\bfseries\scshape}
{\theoremstyle{change}{\theorembodyfont{\normalfont\itshape}
\nthm{conj}[subsection]{Conjecture}  \nthm{*conj}[subsubsection]{Conjecture}
\nthm{**conj}[paragraph]{Conjecture}  \nthm{thm}[subsection]{Theorem}
\nthm{*thm}[subsubsection]{Theorem}  \nthm{**thm}[paragraph]{Theorem}
\nthm{prop}[subsection]{Proposition}  \nthm{*prop}[subsubsection]{Proposition}
\nthm{**prop}[paragraph]{Proposition}  \nthm{lemma}[subsection]{Lemma}
\nthm{*lemma}[subsubsection]{Lemma}  \nthm{**lemma}[paragraph]{Lemma}
\nthm{cor}[subsection]{Corollary}  \nthm{*cor}[subsubsection]{Corollary}
\nthm{**cor}[paragraph]{Corollary}}
{\theorembodyfont{\normalfont\rmfamily}
\nthm{rem}[subsection]{Remark}  \nthm{*rem}[subsubsection]{Remark}
\nthm{**rem}[paragraph]{Remark}  \nthm{defn}[subsection]{Definition}
\nthm{*defn}[subsubsection]{Definition}  \nthm{**defn}[paragraph]{Definition}}}
{\theoremstyle{nonumplain}
{\theorembodyfont{\normalfont\itshape}
\nthm{thm*}[subparagraph]{Theorem}   \nthm{prop*}[subparagraph]{Proposition}
\nthm{lemma*}[subparagraph]{Lemma}   \nthm{cor*}[subparagraph]{Corollary}
\nthm{conj*}[subparagraph]{Conjecture}}
{\theorembodyfont{\normalfont\rmfamily}
 \nthm{rem*}[subparagraph]{Remark}   \nthm{defn*}[subparagraph]{Definition}}}
\DeclareMathAlphabet\eusm{U}{eus}{m}{n}

\def\makeop#1{\expandafter\def\csname#1\endcsname
  {\mathop{\rm #1}\nolimits}\ignorespaces}
\makeop{ad}
\makeop{Ad}
\makeop{Alb}
\makeop{Alt}
\makeop{Aut}
\makeop{Br}
\makeop{can}
\makeop{Card}
\makeop{card}
\makeop{Char}
\makeop{Cl}
\makeop{Coker}
\makeop{coker}
\makeop{Cond}
\makeop{cond}
\makeop{Corr}
\makeop{depth}
\makeop{diag}
\makeop{Disc}
\makeop{disc}
\makeop{Div}
\makeop{div}
\makeop{dom}
\makeop{End}
\makeop{Ext}
\makeop{Fil}
\makeop{Frac}
\makeop{Frob}
\makeop{Fr}
\makeop{Gal}
\makeop{GL}
\makeop{GSpin}
\makeop{GSp}
\makeop{Gr}
\makeop{gr}
\makeop{Hom}
\makeop{hom}
\makeop{Id}
\makeop{id}
\makeop{Im}
\makeop{im}
\makeop{Ind}
\makeop{ind}
\makeop{Int}
\makeop{inv}
\makeop{Irr}
\makeop{irr}
\makeop{Isom}
\makeop{Ker}
\makeop{ker}
\makeop{len}
\makeop{length}
\makeop{Lie}
\makeop{Nr}
\makeop{Nrd}
\makeop{O}
\makeop{Ob}
\makeop{ob}
\makeop{ord}
\makeop{PGL}
\makeop{Pic}
\makeop{Proj}
\makeop{Pr}
\makeop{pr}
\makeop{PSL}
\makeop{rank}
\makeop{RD}
\makeop{Res}
\makeop{Sgn}
\makeop{sgn}
\makeop{SL}
\makeop{SO}
\makeop{Sp}
\makeop{sp}
\makeop{Spec}
\makeop{Spf}
\makeop{Spin}
\makeop{Sp}
\makeop{Stab}
\makeop{stab}
\makeop{std}
\makeop{Std}
\makeop{st}
\makeop{St}
\makeop{Supp}
\makeop{supp}
\makeop{SU}
\makeop{Swan}
\makeop{Sym}
\makeop{Tor}
\makeop{tor}
\makeop{tors}
\makeop{torsion}
\makeop{Trd}
\makeop{Tr}
\makeop{tr}
\makeop{U}
\makeop{val}

\def\makebb#1{\expandafter\def
  \csname bb#1\endcsname{{\mathbb{#1}}}\ignorespaces}
\def\makerm#1{\expandafter\def\csname rm#1\endcsname{{\rm #1}}\ignorespaces}
\def\makebf#1{\expandafter\def\csname bf#1\endcsname{{\bf #1}}\ignorespaces}
\def\makegr#1{\expandafter\def
  \csname gr#1\endcsname{{\mathfrak{#1}}}\ignorespaces}
\def\makescr#1{\expandafter\def
  \csname scr#1\endcsname{{\mathscr{#1}}}\ignorespaces}
\def\makecal#1{\expandafter\def\csname cal#1\endcsname{{\cal #1}}\ignorespaces}
\def\makeudl#1{\expandafter\def\csname udl#1\endcsname{{\underline{#1}}}\ignorespaces}
\def\doLetters#1{#1A #1B #1C #1D #1E #1F #1G #1H #1I #1J #1K #1L #1M
                 #1N #1O #1P #1Q #1R #1S #1T #1U #1V #1W #1X #1Y #1Z}
\def\doletters#1{#1a #1b #1c #1d #1e #1f #1g #1h #1i #1j #1k #1l #1m
                 #1n #1o #1p #1q #1r #1s #1t #1u #1v #1w #1x #1y #1z}
\doLetters\makebb   \doLetters\makecal  \doLetters\makerm  \doLetters\makebf \doLetters\makescr
\doletters\makerm  \doletters\makebf   \doLetters\makegr
\doletters\makegr \doLetters\makeudl  \doletters\makeudl
\def\proof{\medbreak\noindent{\scshape Proof.}\enspace}%
     \def\qed{\qedmark\medbreak}%
    \def\setminus{\smallsetminus}
   
\def\ringO{{\scrO}}      
\def\injto{\hookrightarrow}  
 
\def\surjto{\twoheadrightarrow}

\newcommand{\tworows}[2]{\genfrac{}{}{0pt}{}{#1}{#2}}

\makeatletter   \newdimen\mina@@\mina@@=18pt
\newcommand{\xrtarw}[2][]{\mathrel{\mathop{\,\setbox\z@\vbox{\m@th
  \hbox{$\scriptstyle\;{#1}\;\;$}\hbox{$\m@th\scriptstyle\;{#2}\;\;$}}%
  \hbox to\ifdim\wd\z@>\mina@@\wd\z@\else\mina@@\fi{\rightarrowfill@
  \displaystyle}\,}\limits^{#2}\@ifnotempty{#1}{_{#1}}}}
\newcommand{\xltarw}[2][]{\mathrel{\mathop{\,\setbox\z@\vbox{\m@th
  \hbox{$\scriptstyle\;\;{#1}\;$}\hbox{$\m@th\scriptstyle\;\;{#2}\;\;$}}%
  \hbox to\ifdim\wd\z@>\mina@@\wd\z@\else\mina@@\fi{\leftarrowfill@
  \displaystyle}\,}\limits^{#2}\@ifnotempty{#1}{_{#1}}}}
\makeatother
\def\XYmatrix{\xymatrix@M=5pt} % make \xymatrix not too cluttered
\def\ncmd{\newcommand}
\ncmd{\xysubset}[1][r]{\ar@<-2.5pt>@{^(-}[#1]\ar@<2.5pt>@{_(-}[#1]}
\ncmd{\XYmatrixc}[1]{\vcenter{\XYmatrix{#1}}}
\ncmd{\xyto}[1][r]{\ar@{->}[#1]}      \ncmd{\xyinj}[1][r]{\ar@{^(->}[#1]}
\ncmd{\xysurj}[1][r]{\ar@{->>}[#1]}   \ncmd{\xyline}[1][r]{\ar@{-}[#1]}
\ncmd{\xydotsto}[1][r]{\ar@{.>}[#1]}  \ncmd{\xydots}[1][r]{\ar@{.}[#1]}
\ncmd{\xyleadsto}[1][r]{\ar@{~>}[#1]} \ncmd{\xyeq}[1][r]{\ar@{=}[#1]}
\ncmd{\xyequal}[1][r]{\ar@{=}[#1]}    \ncmd{\xyequals}[1][r]{\ar@{=}[#1]}
\ncmd{\xymapsto}[1][r]{\ar@{|->}[#1]}\ncmd{\xyimplies}[1][r]{\ar@{=>}[#1]}
\ncmd{\xytofrom}[1][r]{\ar@{<->}[#1]} 
  
\def\XYTOTO[#1]^#2_#3{\xyto[#1]<0.5ex>^{#2}\xyto[#1]<-0.5ex>_{#3}}
\ncmd{\xytoto}[1][r]{\XYTOTO[#1]}
\def\beginmat{\begin{pmatrix}}\def\endmat{\end{pmatrix}}
\def\pmat#1]{{\def\beginmat{\begin{pmatrix}}\def\endmat{\end{pmatrix}}\mat#1]}}
\def\bmat#1]{{\def\beginmat{\begin{bmatrix}}\def\endmat{\end{bmatrix}}\mat#1]}}
\def\Bmat#1]{{\def\beginmat{\begin{Bmatrix}}\def\endmat{\end{Bmatrix}}\mat#1]}}
\def\vmat#1]{{\def\beginmat{\begin{vmatrix}}\def\endmat{\end{vmatrix}}\mat#1]}}
\def\Vmat#1]{{\def\beginmat{\begin{Vmatrix}}\def\endmat{\end{Vmatrix}}\mat#1]}}
\def\smat#1]{{\def\beginmat{\begin{smallmatrix}}%
  \def\endmat{\end{smallmatrix}}\left(\mat#1]\right)}}
\def\mat#1#2]{\ifcase#1\or \matA#2]\or \matAA#2]\or \matAAA#2]\fi}
\def\matA  #1#2]{\ifcase#1\or \matAB  #2]\or \matABB  #2]\or \matABBB  #2]\fi}
\def\matAA #1#2]{\ifcase#1\or \matAAB #2]\or \matAABB #2]\or \matAABBB #2]\fi}
\def\matAAA#1#2]{\ifcase#1\or \matAAAB#2]\or \matAAABB#2]\or \matAAABBB#2]\fi}
\def\matAB[#1]{\beginmat#1\endmat}
\def\matABB[#1,#2]{\beginmat#1&#2\endmat}
\def\matABBB[#1,#2,#3]{\beginmat#1&#2&#3\endmat}
\def\matAAB[#1;#2]{\beginmat#1\\#2\endmat}
\def\matAABB[#1,#2;#3,#4]{\beginmat#1&#2\\#3&#4\endmat}
\def\matAABBB[#1,#2,#3;#4,#5,#6]{\beginmat
   #1&#2&#3\\#4&#5&#6\endmat}
\def\matAAAB[#1;#2;#3]{\beginmat#1\\#2\\#3\endmat}
\def\matAAABB[#1,#2;#3,#4;#5,#6]{\beginmat
   #1&#2\\#3&#4\\#5&#6\endmat}
\def\matAAABBB[#1,#2,#3;#4,#5,#6;#7,#8,#9]{\beginmat
   #1&#2&#3\\#4&#5&#6\\#7&#8&#9\endmat}
\def\beginalignorgather#1#2\endalignorgather{\def\tmpx{@}\def\tmpy{#1}%
  \ifx\tmpx\tmpy \begin{align*}#2\end{align*}
  \else\begin{gather*}#1#2\end{gather*}\fi}
\def\[#1\]{\beginalignorgather#1\endalignorgather}

\makeop{Ob}
\makeop{Gr}
\makeop{Ker}
\makeop{ad}

\def\dbltag#1#2{\tag*{\hbox to 0pt{\hbox to \hsize{\hfil
        #2}\hss}#1}}%twosided tag in equations

\newcommand{\lowsim}{\smash{\hbox{\lower2.5pt
  \hbox{\(\scriptstyle\sim\)}}}}%better for putting above an arrow (for isom)
\makeatletter
\def\itemize{%
  \ifnum \@itemdepth >\thr@@\@toodeep\else
    \advance\@itemdepth\@ne
    \edef\@itemitem{labelitem\romannumeral\the\@itemdepth}%
    \expandafter
    \list
      \csname\@itemitem\endcsname
      {\def\makelabel##1{\hss\llap{##1}}}%
  \fi\itemsep=-2pt}
\makeatother
\usepackage[all]{xy}
\def\Quot{\mathop{\rm Quot}}
\def\Ms{\scrM}
\def\Mss{\overline{\scrM}}
\def\isom{\simeq}
\def\tensor{\otimes}
\def\into{\injto}
\def\D{{\mathcal D}}
\def\O{{\ringO}}
\makeop{Diff}
\def\qed{\medbreak}

\begin{document}
\renewcommand{\baselinestretch}{1.00}\normalsize

\begin{center}\scshape
\LARGE
  On vector bundles destabilized by \\Frobenius pull-back\\
\medbreak
\normalsize
\smallbreak
Kirti Joshi, S. Ramanan, Eugene Z. Xia and
  Jiu-Kang Yu\footnote{partially supported by a Sloan Foundation
    Fellowship and by grant DMS 0100678
    from the National Science Foundation}\\
\smallbreak\small
Version 0.99,
July 15, 2002
\end{center}
\bigbreak

\section{Introduction}

Let \(X\) be an irreducible smooth projective curve of genus \(g\) over an
algebraically closed field \(k\) of characteristic \(p>0\), and \(F:X
\to X\) the absolute Frobenius morphism on \(X\).  It is known that
pulling back a stable vector bundle on \(X\) by \(F\) may destroy
stability.  One may measure the failure of (semi-)stability by the
Harder-Narasimhan polygons of vector bundles.

In more formal language, let \(n \geq 2\) be an integer, \(\Ms\)
the coarse moduli space of stable vector bundles of rank \(n\) and
a fixed degree on \(X\). Applying a theorem of Shatz %(an analogue
%of the Grothendieck specialization theorem for \(F\)-isocrystals)
to the pull-back by \(F\) of the universal bundle (assuming the
existence) on \(\Ms\), we see that \(\Ms\) has a canonical
stratification by Harder-Narasimhan polygons (\cite{L}).  We call
this the {\it
  Frobenius stratification.}  This interesting
extra structure on \(\Ms\) is a feature of characteristic \(p>0\).
However, very little is known about the strata of the Frobenius
stratification. Scattered constructions of points outside of the
largest (semi-stable) stratum can be found in \cite{gieseker73},
\cite{RR}, \cite{raynaud82a}. Complete classification of such
points is only known when \(p=2\), \(n=2\), and \(g=2\) by
\cite{joshi00c} and \cite{laszlo00a}.

Our main result here settles the problem for the case of \(p=2\)
and \(n=2\).  On any curve \(X\) of genus \(\geq 2\), we provide a
complete classification of rank-\(2\) semi-stable vector bundles
\(V\) with \(F^*V\) not semi-stable. This also shows that the
bound in \cite[Theorem 3.1]{Su} is sharp. We also obtain fairly
good information about the locus destabilized by Frobenius in the
moduli space, including the irreducibility and the dimension of
each non-empty Frobenius stratum.   In particular we show that the
locus of Frobenius destabilized bundles has dimension $3g-4$ in the
moduli space of semi-stable bundles of rank two. An interesting
consequence of our classification is that high instability of
\(F^*V\) implies high stability of \(V\).

%Moved this part up.
In addition, we show that the Gunning bundle descends
when $g$ is even.  If $g$ is odd, then the Gunning bundle twisted
by any odd degree line
bundle also descends.
%carries such a connection.
%Consolidated the various constructions.
%Our method also \(p=2\)yields less complete results in other cases.

We also construct stable bundles that are destabilized by Frobenius
in the following situations: (1) \(p=2\) and \(n=4\), (2) \(p=n=3\),
(3) \(p=n=5\) and $g \geq 3$.
%we construct stable bundles with
%\(n=3\) that are destabilized by Frobenius provided that $g\geq 5$.
%Similarly, when \(p=5\), we construct such bundles with \(n=5\)
%provided that $g\geq 3$.

%We also construct Frobenius destabilized stable bundles of rank
%four in characteristic two.

The problem studied here can be cast in the generality of
principal \(G\)-bundles over \(X\), where \(G\) is a connected
reductive group over \(k\). More precisely, consider the pull-back
by \(F\) of the universal object on the moduli stack of
semi-stable principal \(G\)-bundles on \(X\). Atiyah-Bott's
generalization of the Harder-Narasimhan filtration should then
give a canonical stratification of the moduli stack (\cite{AB},
see also \cite{C}).

There is a connection between Frobenius destabilized bundles and
(pre)-opers. The investigation of this connection is largely
inspired by \cite{beilinson00b}.
%Our methods in characteristic two
%also show that the Gunning bundle (for even genus) descends under
%carries a connection with $p$-curvature zero; if the genus is odd
%we show that the Gunning bundle twisted by any odd degree line
%bundle carries such a connection.
The new phenomena we
observe here is that, in characteristic $p>0$, pre-opers exist and
they need not be opers (indeed all the examples of pre-opers
provided here are not opers).  When $p=0$, all pre-opers are
opers.

The results of this paper were discovered independently by the first
two authors and the last two authors.  When it was realized that
there was a considerable overlap in the techniques and the
results, we decided to write it jointly. The last two
named authors thank  C.-L. Chai  for discussions and
also wish to thank the hospitality of the National Center for
Theoretical Sciences, Hsinchu, Taiwan.

\section{Generalities}\label{generalities}
\subsection[Notations]\label{notations}
    The following notations are in force throughout this paper
unless otherwise specified.  Let $X$ be a smooth, projective curve of
genus $g\geq 2$ over an algebraically closed field $k$ of characteristic
$p>0$.  Let $\Omega_X^1$ be the sheaf of 1-forms on $X$ and $T_X$
the tangent bundle of $X$.

Let $V$ be a vector bundle on $X$ and denote by $F^*(V)$ the pull-back
of $V$ by $F$. If $V=L$ is a line bundle, then $F^*(L)=L^{\otimes p}$.  We
write $V^*=\Hom_{\O_X}(V,\O_X)$ for the dual bundle of $V$.  Denote by
$\chi(V)$ the Euler
characteristic of $V$.  By Riemann-Roch,
\[
\chi(E)=\deg(V)+\rank(V)(1-g).
\]
Denote by
$\mu(V)=\deg(V)/\rank(V)$,
the {\it slope\/} of $V$.

\subsection[Stability] \label{stablity}
A vector bundle $V$ is {\it stable\/}
(resp.~{\it semi-stable\/}) if for any non-zero sub-bundle $W\subset V$,
$\mu(W)<\mu(V)$ (resp. $\mu(W)\leq \mu(V)$).  A non-zero sub-bundle
\(W\subset V\) with \(\mu(W) \geq \mu(V)\) will be called a {\it
  destabilizing sub-bundle}.
% modified JK, 02/03/2002, property not explicitly used.
%If $W,V$
%are semi-stable and $\mu(W)<\mu(V)$ then
%\[
%\Hom(V,W)=0
%\]
%(\cite[Proposition 5.3.3]{L}).

\subsection[Harder-Narasimhan filtration]\label{harder-narasimhan}
Let $V$ be a vector bundle on $X$. Then there exists a unique
filtration (see \cite[5.4]{L}), called the {\em Harder-Narasimhan
  filtration}, by sub-bundles
\[
0=V_0\subset V_1\subset\cdots \subset V_{h-1}\subset V_h=V
\]
such that $V_i/V_{i-1}$ is semi-stable of slope $\mu_i$ and
\[
\mu_1>\mu_2>\cdots>\mu_h.
\]
The data of \(\bigl(\dim (V_i/V_{i-1}), \mu _i\bigr)\) can be encoded into a
polygon, called the {\it Harder-Narasimhan polygon} (see
\cite[11.1]{L}).  The Harder-Narasimhan polygon can  be regarded
as a measure of instability.

\subsection[A measure of stability]\label{A measure of stability}  Following
\cite{LN}, for a rank-\(2\) vector bundle \(V\), we put
\[
s(V)=\deg(V)-2\max\{\deg(L):L\injto V\},
\]
where the maximum is taken over all rank-\(1\) sub-bundle of
\(V\).  By definition, \(s(V)>0\) (resp.~\(s(V)\geq0\)) if and
only if \(V\) is stable (resp.~semi-stable). When \(s(V)\leq 0\),
the information of \((s(V),\deg(V))\) is the same as that of the
Harder-Narasimhan polygon of \(V\).  Therefore, one may regard
\(s\) as a measure of stability extrapolating the
Harder-Narasimhan polygons, though it is only for the rank-\(2\)
case (for possible variants for the higher rank case, see
\cite{brambila-paz98}; for general reductive group, see
\cite{holla01}).

\subsection{}\label{p-curvature}
For any vector bundle with a connection $(V,\nabla)$, there
exists a $p$-linear morphism of $\O_X$-modules, called the {\em
$p$-curvature of $\nabla$},
\[
\psi:T_X \to \End(V)
\]
which measures the obstruction to the Lie algebra homomorphism
$\nabla:T_X\to \End(V)$ being a homomorphism of $p$-Lie algebras.
A connection is $p$-{\it flat\/} if $\psi $ is zero. A vector bundle is
$p$-{\it flat\/} if it admits a $p$-flat connection.

By a theorem of Cartier (\cite[Theorem 5.1, page 190]{katz70}), there
exists a vector bundle \(W\) on \(X\) such that \(F^*(W)\isom V\) if
and only if $V$ carries a $p$-flat connection
\[ \nabla:V\to\Omega_X^1\tensor V
\]
such that the natural map
\[
F^*(V^{\nabla=0})\to V,
\]
(where $V^{\nabla=0}$ is the module of flat sections
considered as an $\O_X$-module), is an isomorphism.

\subsection{}\label{second-fundamental}
Suppose $(V,\nabla)$ is a vector bundle with a connection and
$W \subset V$ a sub-bundle.  Then there is a natural map (the second
fundamental form)
\[
    T_X\to \Hom(W,V/W),
\]
which is zero if and only if $\nabla$ preserves $W$.  By Cartier's
theorem, if $(V,\nabla)$ is $p$-flat
and $W \into V$ is a sub-bundle preserved by
$\nabla$, then $\nabla$ restricts to a $p$-flat connection on $W$.

\subsection{}\label{raynaud-thm}
Let $B_1$ be the vector bundle defined by the exact sequence
\[
0\to \O_X \to F_*(\O_X)\to B_1\to 0
\]
The bundle $B_1$ is semi-stable of slope $g-1$ (and degree
$(p-1)(g-1)$); moreover, for $p>2$, $F^*(B_1)$ is not semi-stable
\cite{raynaud82a}.  For $p=2$, \(B_1\) is a theta characteristic,
i.e.~\(B_1^{\tensor2}=\Omega^1_X\) \cite{raynaud82a}. By
\cite[Proposition~1.1]{raynaud78}, $B_1$ is stable when $p=3$ and
$g\geq 2$.

\begin{lemma}\label{determinant}
Let $L$ be a line bundle on $X$. Then
\[
\det(F_*L)=\det(B_1)\tensor L.
\]
\end{lemma}
\begin{proof}
See \cite[Chapter 4, Exercise
2.6]{hartshorne-algebraic}.
\end{proof}

\subsection{}\label{mu-calculation}
Let $V$ be a vector bundle on $X$. Then
\[
\deg(F_*V)=\deg(V)+\rank(V)\deg(B_1).
\]
This follows from Riemann-Roch
and the fact that $\chi(F_*V)=\chi(V)$ or by
\ref{determinant}.  In particular
\[
\mu(F_*V)=\frac{1}{p}\mu(V)+(1-\frac{1}{p})(g-1).
\]
\subsection[Duality]\label{duality}
% modified by JK, 02/03/2002
%Let $L\in\Pic(X)$ be a line bundle on $X$.
%Following \cite[section 1.16, page
%70]{ramanathan87}, we have
%\[
%F_*(L)^*\isom F_*(L^{-1}\tensor(\Omega^1_X)^{\tensor(1-p)}).
%\]
%In other words, the dual of $F_*(L)$ is also a push-forward of a line bundle.
Let $V$ be a vector bundle on $X$.
Following \cite[section 1.16, page
70]{ramanathan87}, we have
\[
F_*(V)^*\isom F_*\bigl(V^*\tensor(\Omega^1_X)^{\tensor(1-p)}\bigr).
\]
Thus the dual of $F_*(V)$ is of the form \(F_*(V')\).  We will often
make use of this fact together with the following simple lemma.

\begin{lemma*} Let \(V\) be a vector bundle of rank \(n\) on \(X\),
  \(m\) an integer such that \(0 < m < n\).  The following are
  equivalent:
  \begin{itemize}
  \item [\rm (i)] for all sub-bundle \(W\) of \(V\) of rank \(m\), we
    have \(\mu(W) < \mu(V)\) (resp.~\(\mu(W) \leq \mu(V)\));
  \item [\rm (ii)] for all sub-bundle \(W'\) of \(V^*\) of rank \(n-m\), we
    have \(\mu(W') < \mu(V^*)\) (resp.~\(\mu(W') \leq \mu(V^*)\)).
  \end{itemize}
\end{lemma*}

\section{A general construction}
%modified by K, on 3/27; replaces previous version.
\begin{prop} \label{FFV-is-unstable}
Let $V$ be a vector bundle on $X$.    Then the adjunction map
$F^*(F_*(V))\to V$ is surjective and $\mu(F^*(F_*(V))) = \mu(V) +
(p-1)(g-1)>\mu(V)$.
\end{prop}
\begin{proof}
The surjectivity of the adjunction map is easily check by a local
calculation.  The formula for slope follows from
\ref{mu-calculation}. Hence $\mu(F^*(F_*(V)))
> \mu(V)$.
\end{proof}

% added by Kirti on 3/27
\begin{rem*}
In \ref{canonical-filtration}, we prove a stronger
assertion: $F^*(F_*(V))$ is highly unstable whenever $V$ is
semi-stable.
\end{rem*}

\begin{prop}\label{general-construction}
% notational changes by JK, 02/03/2002
Let $V$ be a semi-stable bundle on $X$.
\begin{itemize}
\item [\rm (i)] For any rank-\(1\) sub-bundle \(L\) of \(F_*V\), we
  have
\[
\mu(L)\leq \mu(F_*V)-\frac{(p-1)(g-1)}{p}.
\]
\item [\rm (ii)] For any rank-\(2\) sub-bundle \(E\) of \(F_*V\), we have
\[
\mu (E)\leq\mu(F_*V)-\frac{1}{p}\left(\frac{pg}{2}-p-g+1\right).
\]
\end{itemize}
\end{prop}

\begin{proof}
If $L\injto F_*V$ is a line sub-bundle, then by adjunction, there
is a non-zero morphism $F^*L\to V$.  Therefore, \(\mu (F^*L)\leq
\mu(V)\); i.e., \(p\cdot \mu (L)\leq\mu(V)=p\cdot \mu
(F_*V)-(p-1)(g-1)\).  Here we have made use of the formulas in
\ref{mu-calculation}.  This proves (i).

Let \(E \injto F_*V\) be a sub-bundle of rank \(2\).  Then by a
theorem of Nagata \cite{holla01}, there is a line sub-bundle \(L
\injto E\) such that \(\mu(L)\geq\mu (E)-g/2\).  Thus we have
\(\mu (E)\leq\mu (L)+g/2\leq \mu(F_*V)-(p-1)(g-1)/p+g/2\).  This
proves (ii).
\end{proof}

%\begin{prop}\label{char3char5}
%{\rm (i)} For any \(p>0\), there is no destabilizing line
%  sub-bundle of \(F_*V\).  {\rm (ii)}
%  Let $V$ be a line bundle. Suppose either \(g\geq 3\) and
%  \(p\geq 5\), or \(g\geq5\) and \(p\geq3\), then there is no
%  destabilizing rank-\(2\) sub-bundle of \(F_* V\).
%\end{prop}

\begin{thm}\label{gen-cons-thm}
{\rm (i)}
Let $p=2$ and
%g \geq 2 is in the general assumption.
% and $g\geq 2$.
$V$ a stable bundle of
rank two and even degree on $X$. Then $F_*(V)$ is a semi-stable bundle
of rank 4 and $F^*(F_*(V))$ is not semi-stable.
{\rm (ii)}~Suppose \(p=3\) (resp.~$g\geq 3$ and $p=5$). Let \(V\) be
a line bundle on \(X\).  Then the
bundle $F_*(V)$ is a stable bundle of rank 3 (resp.~5) and
$F^*(F_*(V))$ is not semi-stable.
\end{thm}

%\begin{cor}\label{general-construction-in-char3and5gen}
%Suppose \(p=3\) and $g \geq 5$ (resp.~$g\geq 3$ and $p=5$). Let \(V\) be
%a line bundle on \(X\).  Then the
%bundle $F_*(V)$ is stable.
%\end{cor}

\begin{proof}
(i) By \ref{duality} and \ref{FFV-is-unstable}, it suffices to show that for any sub-bundle \(E\) of
\(F_*V\) of rank \(\leq 2\), we have \(\mu (E) \leq \mu (F_*V)\).  This
is clear when \(\rank E=1\) by \ref{general-construction} (i).
Suppose that \(\rank E=2\) and \(\mu (E) > \mu (F_*V)\).
The proof of ~\ref{general-construction} (ii) gives a line
bundle \(L \injto E\) such that \(\mu (E)\leq \mu (L)+g/2\leq \mu (F_*V)+1/2\).

The assumption that \(\deg V\) is even implies that \(\mu (F_*V) \in
\frac{1}{2}\bbZ\).  Thus we must have \(\mu (E)=\mu
(L)+g/2=\mu(F_*V)+1/2\).  This gives \(\mu (L)=\frac{1}{2}\mu
(V)\) and \(\mu (F^*L)=\mu (V)\), contradicting the stability of \(V\) as
there is a non-zero morphism \(F^*L \to V\) by adjunction.

% rewritten by JK, 02/03/2002
(ii) By \ref{duality} and \ref{FFV-is-unstable},
it suffices to check that \(F_*V\) does not
have a destabilizing sub-bundle of rank \(\leq 1\) (resp.~\(\leq 2\)).
This is immediate from \ref{general-construction}.
\end{proof}
%By Theorem~\ref{general-construction} we know that $F_*(V)$ does
%not have any devitalizing line sub-bundles. Now by duality for the
%Frobenius morphism, the dual of $F_*(V)^*=F_*(V^*\tensor
%K_X^{-1})$ (we use $p=2$) and so if $F_*(V)$ has a destabilizing
%sub-bundle of rank three then its dual has a devitalizing line
%sub-bundle. Thus we see that if $F_*(V)$ is not semi-stable then
%it is destabilized by a rank two sub-bundle. Assume that $W\to
%F_*(V)$ is such a sub-bundle. So
%$\mu(W)>\mu(F_*(V))=\frac{1}{2}\left(\mu(V)+g-1\right)$. By a
%theorem of Nagata \cite{holla01} $W$ carries a line sub-bundle
%with $\mu(E)\geq \mu(W)-\frac{g}{2}$ and so we see that
%$\mu(F^*(E))\geq \mu(W)-g>2\mu(F_*(V))-g=\mu(V)+g-1-g$. As
%$\mu(E)$ is an integer (as $E$ is a line bundle) we see that
%$\mu(F^*(E))>\mu(V)-1$, gives $\mu(F^*(E))\geq \mu(V)$ (here we
%use the fact that $\mu(V)$ is an integer) which contradicts the
%stability of $V$. This completes the proof.
%\end{proof}

\section{A detailed study of the case of rank \(2\) and characteristic
  \(2\)}\label{detailed-char2}

Throughout this section, \(p=2\).  We present our main results on the
classification of rank-\(2\) vector bundles destabilized by Frobenius,
as well as the geometry of the Frobenius stratification.

\subsection[A result on the Gunning bundle]  We begin with an
interesting observation about Gunning extensions, though this
result is not needed in the sequel. Recall that \(B_1\) is a
theta-characteristic \cite[\S4]{raynaud82a}. The unique
non-trivial extension \(0 \to B_1 \to W \to B_1^{-1} \to 0\) is
called {\it the Gunning extension} and the bundle \(W\) is called {\it
  the Gunning bundle.}

\begin{prop*}  Let \(\xi\) be a line bundle and \(V=F_*(\xi\otimes B_1^{-1})\).
The extension
\[\tag{*}
0 \to \xi\otimes B_1 \to F^*V \to \xi\otimes B_1^{-1} \to 0
\]
defines a class in \(\Ext^1(\xi\otimes B_1^{-1}, \xi\otimes
B_1)\simeq H^1(X,B_1^2)\simeq k\).  This class is trivial precisely
when \(\deg(\xi\otimes B_1^{-1})\) is even.
\end{prop*}

\proof Suppose that \(\deg(\xi\otimes B_1^{-1})\) is even.  Then we can
write \(L=\xi\otimes B_1^{-1}=M^2\).  By \cite[\S2]{joshi00c}, there is an
exact sequence \(0 \to M \to V \to M \otimes B_1 \to 0\).  Pulling
back by \(F\), we get \(0 \to L \to F^*V \to L \otimes B_1^2 \to
0\).  This shows that \((*)\) is split.

Suppose that \(L=\xi\otimes B_1^{-1}\) has odd degree \(2n+1\).  By a
theorem of Nagata (\cite{LN}, Cf.~Remark \ref{moduli}), there is an
exact sequence \(0 \to M_1 \to V \to M_2\to0\), where \(M_1,M_2\) are
line bundles with degrees \(n\) and \(n+g\) respectively.  From the
exact sequence \(0 \to M_1^2 \to F^*V \to M_2^2\to 0\), we deduce that
\(\dim\Hom(L,F^*V)\leq \dim\Hom(L, M_1^2) +\dim\Hom(L, M_2^2) =0+g=g\)
by the Riemann-Roch formula.  Since \(\Hom(L,\xi\otimes
B_1)=H^0(X,B_1^2)\) has dimension \(g\), any morphism \(L \to F^*V\)
factors through the sub-module \(\xi\otimes B_1\) in \((*)\).
Therefore, \((*)\) is not split.\qed

\begin{cor*}\label{descent-thm} Let \(W\) be the Gunning bundle and
  \(\xi\) a line bundle of degree \(\equiv g\pmod{2}\).
Then there exists a stable bundle $V$ such that \(F^*V\simeq W \otimes
\xi\). In particular, if \(g\) is even, then the Gunning bundle \(W\)
is the Frobenius pull-back of a stable bundle.
\end{cor*}

% rewritten by JK, 02/03/2002
%\begin{cor*}\label{descent-thm}  Suppose that \(g\) is even.
%  Then there exists a stable bundle $V$ of trivial determinant, whose
%  Frobenius pull-back is the Gunning extension. If $g\geq 2$ is
%  odd then for any line bundle $\xi$ of odd degree, the twisted
%  Gunning bundle $W\tensor\xi$ is the Frobenius pull-back of a
%  stable bundle.
%\end{cor*}
%\begin{proof}
%If $g$ is even we put \(\xi=\ringO_X\) in the above proposition.
%If $g$ is odd then for any line bundle $\xi$ of odd degree
%$V=F_*(\xi\tensor B_1^{-1})$ is stable of determinant $\xi$ and
%$F^*(V)\tensor \xi^{-1}\simeq W$.
%\end{proof}

\begin{rem*} This corollary is implicit in \cite{laszlo00a} in
 the case of an ordinary curve with $g=2,p=2$.
In \cite{gieseker73}, Gieseker proved (by different methods)
an analogous result in any characteristic when $X$ is a
Mumford curve.
\end{rem*}

\subsection[The basic construction]\label{cons}
Henceforth, fix an integer \(d\).
For an injection \(V' \injto
V''\) of vector bundles of the same rank, define the {\it
  co-length\/} \(l\) of \(V'\) in \(V''\) to be the length of the torsion
\(\ringO_X\)-module \(V''/V'\).  Clearly, \(s(V')\geq s(V'')-l\).

We now give a basic construction of stable vector bundles
\(V\) of rank \(2\) with \(F^*V\) not semi-stable.  Let \(l \leq g-2\) be a
non-negative integer, \(L\) a line bundle of degree \(d-1-(g-2-l)\),
and \(V\) a sub-module of \(F_*L\) of co-length \(l\), then
\(\deg V=d\) and \(s(V)\geq(g-1)-l>0\) by \ref{general-construction}.
Therefore, \(V\) is stable.

On the other hand, by adjunction, there is a morphism \(F^*V \to L\), and
the kernel is a line bundle of degree \(\geq
d+1+(g-2-l)>d=\deg(F^*V)/2\).  Therefore,
\(F^*V\) is not semi-stable.

\subsection[Exhaustion] \label{exh}
Suppose that \(V\) is
semi-stable of rank \(2\) and \(F^*V\) is not semi-stable.

Let \(\xi=\det(V)\) and \(d=\deg \xi=\deg V\).  Since \(F^*V\) is
not semi-stable and of degree \(2d\), there are line bundles \(L,L'\) and
an exact sequence \(0 \to L' \to F^*V \to L \to 0\)
with \(\deg L'\geq d+1\), \(\deg L \leq d-1\). By adjunction, this provides a non-zero
morphism \(V \to F_*L\).  If the image is a line
bundle \(M\), we have \(\deg M\geq d/2\) by semi-stability of \(V\),
and \(\deg M \leq (d-1+g-1)/2-(g-1)/2=(d-1)/2\) by
\ref{general-construction}.  This is a contradiction.

Thus the image is of rank \(2\).  Since \(\deg V=d\) and
\(\deg (F_*L)\leq d+(g-2)\), \(V\) is a sub-module
of \(F_*L\) of co-length \(l \leq g-2\), and \(\deg
L=d-1-(g-2-l)\).

Thus the basic construction yields all semi-stable
vector bundles \(V\) of rank
\(2\), with \(F^*V\) not semi-stable.

\begin{*cor} \label{exh1}
If \(V\) is semi-stable of rank \(2\) with \(F^*V\) not semi-stable,
then \(V\) is actually stable.\qed
\end{*cor}

\begin{*cor} \label{exh2}
  The basic construction with \(l=g-2\) already yields all
  semi-stable vector bundles \(V\) of rank \(2\), with \(F^*V\) not semi-stable.
\end{*cor}

\proof
In fact, if \(l<l'\leq g-2\) and \(L'=L\otimes \ringO(D)\) for some
effective divisor \(D\) of degree \(l'-l\) on \(X\),
then \(V \injto F_*L \injto F_* L'\).  Hence \(V\) is also
a sub-module of \(F_*L'\) of co-length \(l'\).  Thus \(V\) arises from
the basic construction with \((l',L')\) playing the role of \((l,L)\).\qed

\subsection[Classification]\label{cls}
Let \(L\) be a line bundle and let
\(Q=Q_l=Q_{l,L}=\Quot_l(F_*L/X/k)\) be the scheme classifying
sub-modules of \(F_*L\) of co-length \(l\) (\cite[3.2]{FGA}).  Let
\[
\scrV \injto\ringO_Q \boxtimes F_*L =(\id \times F)_*(\ringO_Q \boxtimes L)
\]
(sheaves on \(Q \times X\)) be the universal object on \(Q\).  By
adjunction, we have a morphism \((1\times F)^*\scrV \to
\ringO_Q\boxtimes L\).  Let \(\scrF\) be the cokernel.  Then
\(\pr_*\scrF\) is a coherent sheaf on \(Q\), where \(\pr:Q\times X \to
Q\) is the projection (\cite[II.5.20]{hartshorne-algebraic}).  By
\cite[III.12.7.2]{hartshorne-algebraic}, the subset
\[
\{q \in Q : \dim_{\kappa (q)}((\pr_*\scrF) \otimes \kappa(q)
> 0\}
\] is closed.  Its complement is an open sub-scheme, denoted by
\(Q^*=Q_l^*=Q_{l,L}^*\), of \(Q\).  Then \(Q^*\) parameterizes those
\(V\)'s with surjective \(F^*V \to L\).

Let \(\Mss\) be the coarse moduli space of rank-\(2\) semi-stable vector
bundles of degree \(d\) on \(X\).  Let \(\Ms\) be the open sub-scheme
parameterizing stable vector bundles and and \(\Ms_1(k)\subset
\Mss(k)\) the subset of those \(V\)'s such that \(F^*V\) is not
semi-stable.  By \ref{exh1}, \(\Ms_1(k) \subset \Ms(k)\).

\begin{prop*}
The basic construction gives a bijection
\[
\coprod_{\tworows{0\leq l\leq g-2}{\deg L=d-1-(g-2-l)}} Q_{l,L}^*(k) \to
  \Ms_1(k),
\]
where the disjoint union is taken over all \(l\in [0,g-2]\) and a set
of representatives of all
isomorphism classes of line bundles \(L\) of degree \(d-1-(g-2-l)\).
\end{prop*}

\proof By \ref{exh}, the map is a surjection.  Now suppose that
\((l,L,V\subset F_*L)\) and \((l',L',V'\subset F_* L')\) give the same
point in \(\Ms_1(k)\), i.e.~\(V \simeq V'\).  Since the unstable
bundle \(F^*V\) has a unique quotient line bundle of degree
\(<\deg(V)/2\) (i.e.~the second graded piece of the Harder-Narasimhan
filtration), which is isomorphic to \(L\), we must have \(L=L'\).
Consider the diagram
\[
\XYmatrix{
F^* V \xyto \xyto[d]_{\wr} &L\xyequal[d]\\
F^* V' \xyto &L',}
\]
where the vertical arrow is induced from an isomorphism \(V
\xrtarw{\sim} V'\) and the horizontal arrows are the unique quotient
maps.  This diagram is commutative up to a multiplicative scalar in
\(k^*\).  By adjunction, \(V \injto F_*L\) and \(V' \injto F_* L\) have
the same image.  In other words, \(V = V'\) as sub-modules of
\(F_*L\).  This proves the injectivity of the map.\qed

\subsection[Frobenius Stratification] \label{moduli}
To ease the notation, let
\(d_l=d-1-(g-2-l)\).  Let \(\Pic^{d_l}X\) be the moduli space of line
bundles of degree \(d_l\) on \(X\), and \(\scrL \to \Pic^{d_l}(X)
\times X\) the universal line bundle.

By \cite[3.2]{FGA}, there is a scheme \(\scrQ=\scrQ_l=\Quot_l
\bigl((\id \times F)_*\scrL/ (\Pic^{d_l}(X) \times X)/\Pic^{d_l}
X\bigr) \xrtarw{\pi } \Pic^{d_l} X\) such that \(\scrQ_x\) (the fiber
at \(x\)) is \(Q_{\scrL_x}\) for all \(x \in (\Pic^{d_l}X)(k)\).  By
the same argument as before, there is an open sub-scheme
\(\scrQ^*\subset \scrQ\) such that \(\scrQ^*_x=Q^*_{\scrL_x}\) for all
\(x \in \Pic^{d_l}(X)(k)\).  The scheme \(\scrQ \) is projective over
\(\Pic^{d_l}(X)\) (\cite[3.2]{FGA}), hence is proper over \(k\).  By
checking the condition of formal smoothness (cf.\ \cite[8.2.1]{L}),
it can be shown that \(\scrQ \) is smooth over \(\Pic^{d_l}(X)\),
hence is smooth over \(k\).

The Frobenius stratification on the coarse moduli scheme \(\Ms\) is
defined canonically using  Harder-Narasimhan polygons of Frobenius
pull-backs . Concretely, for \(j\geq0\), let \(P_j\) be the
polygon from \((0,0)\) to \((1,d+j)\) to \((0,2d)\).  Let
\(\Ms_0=\Mss\), and for \(j\geq1\), let \(\Ms_j(k)\) be the subset
of \(\Mss(k)\) parameterizing those \(V\)'s such that the
Harder-Narasimhan polygons (\cite[11.1]{L}) of \(F^*V\) lie above
or are equal to \(P_j\). Notice that \(\Ms_1(k)\) agrees with the
one defined in \ref{cls}.

As mentioned in the introduction, the existence of a universal bundle
on \(\Ms\) would imply that each \(\Ms_j(k)\) is Zariski closed by
Shatz's theorem \cite[11.1, last remark]{L}.  In general, one can show
that \(\Ms_j(k)\) is closed by examining the GIT (geometric invariant
theory) construction of \(\Mss\).  This fact also follows from our
basic construction:

\begin{thm*} The subset \(\Ms_j(k)\) is Zariski closed in \(\Mss(k)\), hence underlies a
  reduced closed sub-scheme \(\Ms_j\) of \(\Mss\).
  The scheme \(\Ms_j\) is proper.  The Frobenius
  stratum \(\Ms_j \setminus \Ms_{j+1}\) is non-empty precisely when \(0
  \leq j \leq g-1\).  For \(1 \leq j \leq g-1\), write \(l=g-1-j\).
  Then there is a canonical morphism
\[
\scrQ_l \to \Mss
\]
which has scheme-theoretic image \(\Ms_j\) and induces a bijection
from \(\scrQ_l^*(k)\) to \(\Ms_j(k) \setminus \Ms_{j+1}(k)\).
\end{thm*}

\proof Suppose \(0 \leq l \leq g-2\) and \(j+l=g-1\).  The universal
object \(\scrV \to \scrQ_l \times X\) is a family of stable vector
bundles on \(X\).  This induces a canonical morphism \(\scrQ_l \to
\Mss\).  The image of \(\scrQ_l(k)\) is precisely \(\Ms_j(k)\) by (the
proof of) \ref{exh2}.  Since \(\scrQ_l\) is proper,
\(\Ms_j\) is proper and closed in \(\Mss\).  The rest of the
proposition follows from \ref{cls} and \ref{exh}, and the fact
that \(\scrQ_l^*(k)\) is non-empty for \(0 \leq l \leq g-2\) (see
\ref{dense}).\qed

\begin{rem*}
  By a theorem of Nagata (\cite{LN}, \cite{holla01}), \(s(V)\leq g\)
  for all \(V\).  Therefore, \(s(V)\leq g\) if \(\deg V \equiv g
  \pmod{2}\), and \(s(V)\leq g-1\) if \(\deg V \not\equiv g
  \pmod{2}\).  By \ref{general-construction},
  \(V=F_*L\) achieves the maximum value of \(s\) among rank-\(2\)
  vector bundles of the same degree.

  By the preceding theorem, vector bundles of the form \(V=F_*L\) are
  precisely members of the smallest non-empty Harder-Narasimhan
  stratum \(\Ms_{g-1}\).  Therefore, in a sense \(V\) is most stable
  yet \(F^*V\) is most unstable.  More generally, for \(1 \leq j \leq
  g-1\), we have (from \ref{cons})
\[
s(\Ms_j(k)) \geq \begin{cases}
j &\mbox{if \(d\equiv j \pmod{2}\),}\\
j+1 &\mbox{if \(d\not\equiv j\pmod{2}\).}
\end{cases}
\]
Therefore, high instability of \(F^*V\) implies high stability of
\(V\).
\end{rem*}

\subsection[Irreducibility]\label{irr}
We will make use of the following simple lemma.

\begin{*lemma} \label{irr-cri}
Let \(Y\) be a proper scheme over \(k\), \(S\) an
irreducible scheme of finite type over \(k\) of dimension \(s\), \(r\) an integer
\(\geq 0\), and \(f: Y \to S\)
a surjective morphism.  Suppose that all fibers of \(f\) are
irreducible of dimension \(r\).  Then \(Y\) is irreducible of
dimension \(s+r\).\qed
\end{*lemma}

\begin{*lemma}\label{irrQ}
The scheme \(\scrQ=\scrQ_l\) is irreducible of dimension \(2l+g\).
\end{*lemma}

\proof
There is a surjective morphism
(\cite[\S6]{FGA})
\[
\delta :\scrQ \to
\Div^l(X)=\Sym^l(X), \qquad q\mapsto\sum_{P \in X(k)}
\length_{\ringO_P}\bigl((F_*\scrL_{\pi(q)})/\scrV_q\bigr) \cdot P.
\]
The morphism \(\scrQ \to \Div^l(X) \times \Pic^{d_l}(X)\) is again a
surjection.  The fibers are irreducible schemes of dimension \(l\)
according to the last lemma of \cite{markman02}.  Since \(\scrQ\) is proper,
the result follows from \ref{irr-cri}. \qed

\begin{*lemma}\label{dense}
 \(\scrQ^*\) is open and dense in \(\scrQ\).
\end{*lemma}

\proof By the construction in \ref{cls} and \ref{moduli},
\(\scrQ^*\) is open in \(\scrQ\).  Since \(\scrQ\) is irreducible of
dimension \(2l+g\), it suffices to show that \(\scrQ^*\) is non-empty.
We will do more by exhibiting an open subset of \(\scrQ^*\) of dimension \(2l+g\).

Indeed, let \(B(X,l) \subset \Div^l(X)\) be the open sub-scheme parameterizing
multiplicity-free divisors of degree \(l\), also known as the
configuration space of unordered \(l\) points in \(X\).  Let
\(U\) be the inverse image of \(B(X,l) \times \Pic^{d_l}(X)\) under
\(\scrQ^* \to \Div^l(X) \times \Pic^{d_l}(X)\).  A quick
calculation shows that each fiber of \(U \to B(X,l) \times \Pic^{d_l}(X)\) is
isomorphic to \(\bbA^l\).  Therefore, \(U\) is an open subset of
\(\scrQ^*\) of dimension \(2l+g\). \qed

\begin{*thm}
  For \(1 \leq j \leq g-1\), \(\Ms_j\) is proper, irreducible, and of dimension
  \(g+2(g-1-j)\).  In particular, \(\Ms_1\) is irreducible and of
  dimension \(3g-4\).\qed
\end{*thm}

\subsection[Fixing the determinant]
Fix a line bundle \(\xi\) of degree \(d\).
Let \(\Mss(\xi )\subset \Mss\) be the closed sub-scheme of \(\Mss \) parameterizing
those \(V\)'s with \(\det(V)=\xi\).  Let \(\Ms
_j(\xi)=\Mss(\xi)\cap \Ms_j\) for \(j\geq0\).

\begin{rem*}
For \(1 \leq j \leq g-1\), \(\dim \Ms_j(\xi)=2(g-1-j)\).  In particular, \(\dim
\Ms_1(\xi)=2(g-2)\).
\end{rem*}

\proof Since \(\Ms_j(\xi)\) is nothing but the fiber of the surjective
morphism \(\det: \Ms_j \to \Pic^d(X)\), it has dimension \(2(g-1-j)\)
for a dense open set of \(\xi \in \Pic^d(X)(k)\).  However,
\(\Ms_j(\xi _1)\) is isomorphic to \(\Ms_j(\xi_2)\) for all \(\xi_1,\xi_2\in
\Pic^d(X)(k)\), via \(V \mapsto V \otimes L\), where \(L^2\simeq
\xi_2\otimes\xi_1^{-1}\).  Thus the remark is clear.\qed

A slight variation of the above argument shows that \(\Ms_j(\xi)\)
is irreducible.  Alternatively, assume \(1 \leq j \leq g-1\).  Let
\(l=g-1-j\) and let \(\scrQ(\xi)=\scrQ_l(\xi)\) be the inverse image
of \(\xi\) under \(\scrQ \to \Pic^d(X)\), \(q \mapsto \det(\scrV_q)\).
Since \(\det(\scrV_q)=B_1 \otimes \scrL_{\pi(q)} \otimes
\ringO(-\delta(q))\), the morphism \(\det: \scrQ \to \Pic^d(X)\)
factors as
\[
\scrQ \to \Div^l(X) \times \Pic^{d_l}(X) \xrtarw{\psi} \Pic^d(X),
\]
where \(\psi\) is \((D,L) \mapsto B_1 \otimes L \otimes \ringO(-D)\).
It is clear that \(\psi^{-1}(\xi)\) is isomorphic to \(\Div^l(X)\),
and hence is an irreducible variety.

The fibers of \(\scrQ(\xi) \to \psi^{-1}(\xi)\) are just some
fibers of \(\scrQ \to \Div^l(X) \times \Pic^{d_l}(X)\); hence they are
irreducible of dimension \(l\) as in the proof of
\ref{irrQ}. Being a closed sub-scheme of \(\scrQ\),
\(\scrQ(\xi)\) is proper, thus,
irreducible by \ref{irr-cri}.  Now it is easy to deduce

%the meaning of ``strata'' in the previous version does not conform to
%the usual definition.
\begin{thm*} There is a canonical (Frobenius) stratification by
Harder-Narasimhan polygons
%on the scheme \(\Mss(\xi)\)
%has strata
\[
\emptyset=\Ms_{g}(\xi)\subset \Ms_{g-1}(\xi) \subset \cdots
\subset \Ms_0(\xi)=\Mss (\xi)
,
\]
with \(\Ms_j(\xi)\) non-empty, proper, irreducible, and of
dimension \(2(g-1-j)\) for \(1 \leq j \leq g-1\).\qed
\end{thm*}

\subsection[A variant] Let \(\Ms '(k)\) be the subset of \(\Mss(k)\)
consisting of those \(V\) such that \(F^*V\) is not stable.  Clearly,
\(\Ms '(k) \supset \Ms_1(k)\).

By \ref{exh1}, the closed subset \(\Ms^{\rm
  ns}(k)=\Mss(k)\setminus \Ms(k)\) is contained in \(\Ms
'(k)\setminus \Ms_1(k)\).  On the other hand, if \(V \in \Ms
'(k)\setminus \Ms^{\rm ns}(k)\), the argument of \ref{exh} shows
that there is a line bundle \(L\) of degree \(d\) such that \(V \injto
F_*L\) is a sub-module of co-length \(\leq g-1\).  Conversely, the
argument of \ref{cons} shows that if \(V\) is of co-length \(\leq
g-1\) in \(F_*L\) for some \(L\) of degree \(d\), then \(V \in \Ms
'(k)\).

Thus we conclude that \(\Ms '(k)\) is the union of \(\Ms^{\rm
  ns}(k)\) and the image \(\Ms'_0(k)\) of \(\scrQ_{g-1}(k)\) for a
suitable morphism \(\scrQ_{g-1} \to \Mss\), where \(\scrQ_{g-1}\) is defined in
\ref{moduli}.  It follows that \(\Ms'_0(k)\) and \(\Ms '(k)\)
are Zariski closed in \(\Mss(k)\), hence are sets of \(k\)-points
of reduced closed sub-scheme \(\Ms '_0\) and \(\Ms '\) of
\(\Mss\).

\begin{thm*} The scheme \(\Ms '_0\) is irreducible of dimension \(3g-2\).  It
  contains two disjoint closed subsets: \(\Ms '_0 \cap \Ms^{\rm
    ns}\), which is irreducible of dimension \(2g-1\) when \(d\) is even
  and empty when \(d\) is odd, and \(\Ms_1\), which is irreducible of
  dimension \(3g-4\).
\end{thm*}

\begin{rem*} \(\Ms ' \setminus \Ms_1\) is the first stratum in the
  \(s\)-stratification (\cite{LN}) which is not a Harder-Narasimhan
  stratum.  The other \(s\)-strata are more complicated and not
  pursued here.
\end{rem*}

\proof Since \(\scrQ_{g-1}\) is irreducible, \(\Ms '_0\) is
irreducible.  We now analyze \(\Ms '_0\cap \Ms^{\rm ns}\).
Suppose that \(V \in \Ms '_0(k) \cap \Ms^{\rm ns}(k)\).  Then
\(d=\deg V\) is even and there exists \(L\) of degree \(d\) such that \(V\)
is a sub-module of \(F_*L\) of co-length \(g-1\).  By assumption,
there is a sub-bundle \(M\) of \(V\) of degree \(d/2\).  Adjunction
applied to the composition \(M \injto V \injto F_*L\) provides
a non-zero morphism \(F^*(M)=M^2 \to L\).  This implies that \(M^2 \simeq L\).
We may assume that \(L=M^2\).  Since there is only one (modulo
\(k^*\)) non-zero morphism \(M^2 \to L\), there is only one non-zero
morphism \(M \to F_*(F^* M)\).  By \cite[\S2]{joshi00c}, this morphism is
part of an exact sequence \(0 \to M \to F_*(F^*M) \to M \otimes B_1 \to
0\).   Thus to have \(V\) is to have a sub-module of \(M \otimes B_1\)
of co-length \(g-1\).  Conversely, starting with a sub-module of \(M \otimes
B_1\) of co-length \(g-1\), we obtain a vector bundle \(V \in \Ms
'_0(k) \cap \Ms^{\rm ns}(k)\) as the inverse image of that sub-module
in \(F_*(F^*M)\).

The sub-modules of \(M \otimes B_1\) of co-length \(g-1\) are of the
form \(M \otimes B_1 \otimes \ringO(-D)\) for \(D \in
\Div^{g-1}(X)(k)\).  Thus there is a morphism
\(\pi ': \scrQ '=\Div^{g-1}(X) \times \Pic^{d/2}(X) \to \Mss\)
inducing a surjection \(\scrQ '(k) \to \Ms '_0(k) \cap \Ms^{\rm
  ns}(k)\).  We claim that this morphism is generically finite of
separable degree at most \(2\).  This claim implies that \(\Ms '_0
\cap \Ms^{\rm ns}\) is irreducible of dimension \(2g-1\).

Indeed, there is an open subset \(U\) of \(\Div^{g-1}(X)(k)\) such
that if \(D, D' \in U\) are distinct, then \(D \not\sim D'\).  We now
show that \(\pi '|(U\times \Pic^{d/2}(X)(k))\) is at most 2-to-1.
Suppose that \(D \in U\), \(M \in \Pic^{d/2}(X)(k)\), and \(\pi
'(D,M)=V\).  Then \(V\) has at most two isomorphism classes of
rank-\(1\) sub-bundles of degree \(d/2\), and \(M\) is one of them.
After obtaining \(M\), one can determine \(D\) uniquely by the
condition \(\det(V)\simeq M^2\otimes B_1\otimes \ringO(-D)\).  This
proves the claim.

Next, we consider the morphism \(\scrQ_{g-1} \to
\Ms_0'\).  It induces a surjection \(\scrQ_{g-1}^*(k) \surjto
\Ms_0'(k)\setminus \Ms_1(k)\).  Again the claim is that the morphism
is generically finite of separable degree at most \(2\).
This claim implies that \(\Ms_0'\) is irreducible of dimension
\(3g-2\).

Indeed, let \(U\) be the open subset of \(\scrQ_{g-1}^*(k)\) consisting
of those \(q\)'s such that
\(\ringO(2\delta(q))\not\simeq\Omega_{X/k}\).  Now assume that \(q\in U\) gives
rise to \(V \in \Ms_0'(k)\).  Then there is an exact sequence \(0 \to
L \otimes B_1^2 \otimes \ringO(-2\delta(q)) \to F^*V \to L \to 0\), where
\(L=\scrL_{\pi(q)}\).  The assumption on \(q\) implies that \(F^*V\)
has at most \(2\) quotient line bundles of degree \(d\), say \(F^*V \to L_1\) and
\(F^*V\to L_2\).  Then \(q\) must be one of the two data \(V \injto
F_*L_1\) or \(V \injto F_* L_2\) provided by adjunction.  This proves
the claim.\qed

\subsection[Example] When \(g=2\), \(\Ms _1(\xi)\) is a single point,
corresponding to the vector bundle \(F_*(\xi \otimes B_1^{-1})\).

When \(\xi=B_1\), this refines a result of
\cite[1.1]{joshi00c}, which says that \(\Ms_1(\xi)\) is a single
\(\Pic(X)[2]\)-orbit.

When \(\xi=\ringO_X\), this extends a theorem of Mehta \cite[3.2]{joshi00c},
which states that there are only finitely many rank-\(2\) semi-stable
vector bundles \(V\)'s on \(X\) with \(\det(V)=\ringO_X\) and \(F^*V\)
not semi-stable when \(p\geq3,g=2\).  We now have this result for
\(p=2,g=2\) with the stronger conclusion of uniqueness.

%\subsection[Erratum for \cite{joshi00c}]
%We correct a minor error in the statement of \cite[Theorem
%1.1]{joshi00c}.  The expression ``\(V_1 \in \Ext^1(L_\theta,\ringO_X)\)''
%should be replaced by ``\(V_1 \in S_\theta\)'' (the original version
%is valid when \(L_\theta=B_1\)).  Also, \(\Omega\) should be replaced
%by \(L_\theta\).

\section{Pre-opers and opers}\label{pre-opers-and-opers}
    This section is largely inspired by the work Beilinson and
Drinfel'd \cite{beilinson00b}.
%Pre-opers are a
%characteristic $p>0$ phenomenon which arises from the fact that the
%ring of differential operators on a variety is not generated by
%first order operators.
We show that
pre-opers with connections of $p$-curvature zero provide, under
additional assumptions, examples of Frobenius destabilized bundles.
In small characteristics we describe the lowest Frobenius stratum in
terms of pre-opers.
%A general reference for the formalism of
%opers is \cite[Section 3, page 61]{beilinson00b}.

\subsection[Pre-opers]\label{pre-opers}
Let $V$ be a vector
bundle on $X$ with a flat connection $\nabla$. Suppose that
$\{V_i\}_{0\leq i\leq l}\subset V$ is an increasing filtration by sub-bundles such
that
\begin{enumerate}
\item \(V_0=0\), \(V_l=V\),
\item $\nabla(V_i)\subset V_{i+1}\tensor \Omega^1_X$ for \(0 \leq i
  \leq l-1\),
\item $V_i/V_{i-1}\xrtarw{\nabla} (V_{i+1}/V_i)\tensor
\Omega^1_X$ is an isomorphism for \(1 \leq i\leq l-1\).
\end{enumerate}
Then $(V,\nabla,\{V_i\})$ is said to be a {\it pre-oper}.
A pre-oper is $p$-{\it flat\/} if $\nabla$ has
$p$-curvature zero.

\begin{rem*}
  Let $(V,\nabla,\{V_i\}_{0\leq i\leq l})$ be a pre-oper.  If \(g\geq2\) and \(V_1/V_0\) is
  semi-stable, then the filtration $\{V_i\}_{0\leq i\leq l}$ is nothing but
  the Harder-Narasimhan filtration of \(V\).
\end{rem*}

%\begin{*rem}
%defined for any reductive group $G$. Likewise
%(Note that one can likewise define
%pre-opers for any reductive group $G$ \cite[Section 3, page
%61]{beilinson00b}).
% In this paper
%will not consider pre-opers (and opers) for any other group, so we
%have suppressed the group from our discussions.
%\end{*rem}

%\begin{*rem}
%    Here is a proto-type of a pre-oper in any characteristic
%$p>0$: the bundle $F^*(F_*(\O_X))$ can be identified with the
%space of jets of order at most $p-1$ (see \cite{raynaud82a}) and
%the filtration by order of the differential operator is the
%Harder-Narasimhan filtration.  Then $F^*(F_*(\O_X))$ is a
%pre-oper with the Cartier connection
%and the Harder-Narasimhan filtration.
%\end{*rem}

\subsection[Opers]\label{opers}
%    Let $X$ be a smooth, projective curve over a perfect field of
%characteristic $p>0$.
Let $\D_X=\Diff(\O_X,\O_X)$ be the ring of
differential operators on $X$.
%We caution the reader that $\D_X$
%is not generated by first order operators in characteristic $p>0$.
An {\it oper} $(V,\nabla,\{V_i\}_{0\leq i\leq l})$ is a pre-oper such
that the connection $\nabla$ on $V$ extends to a structure of
$\D_X$-module on $V$.

\begin{rem*} (i)
    By a Theorem of Katz (see \cite[Theorem~1.3, page
4]{gieseker75}), $V$ is $\D_X$-module if and only if there exists a
sequence of vector bundles $V^{i}$ such that $F^*(V^1)=V$ and
$F^*(V^{i+1})=V^i$. In particular, if $(V,\nabla)$ is a vector
bundle such that $\nabla$ extends to a $\D_X$-module structure on
$V$ then $\nabla$ has $p$-curvature zero, so any oper is
automatically $p$-flat.

(ii) When \(l=\rank V\), what we called an oper here is the same as
an $GL_l$-oper as defined in \cite{beilinson00b}.
\end{rem*}

\subsection[A canonical filtration]\label{canonical-filtration} Let \(W\) be a vector bundle on
\(X\).  We define a canonical increasing filtration on \(V=F^*(F_*W)\)
by abelian sub-sheaves \(\{V_i\}_{0\leq i\leq p}\), as follows:
\[@
V_p&=V,\\
V_{p-1}&=\ker(V_p =F^*(F_*W) \to W),\\
V_i&=\ker\bigl(V_{i+1} \xrtarw{\nabla^{\rm Cartier}} V\otimes\Omega_X^1\to
(V/V_{i+1})\otimes\Omega_X^1\bigr),\qquad 0\leq i\leq p-2.
\]
It is elementary to check by induction that each \(V_i\) is actually
an \(\ringO_X\)-sub-module of \(V\).

\begin{thm*}
\begin{itemize}
  \item [\rm (i)] \(V_p/V_{p-1}\) is isomorphic to \(W\).
  \item [\rm (ii)] \(\bigl(V,\nabla^{\rm Cartier},\{V_i\}_{0\leq i\leq p}\bigr)\)
    is a pre-oper.
  \item [\rm (iii)] If \(g \geq 2\) and \(W\) is semi-stable, \(\{V_i\}_{0\leq i\leq p}\) is simply the
    Harder-Narasimhan filtration on \(V\).
\end{itemize}
\end{thm*}

\proof It is clear that (i) and (ii) imply (iii).  To prove (i) and (ii), we notice
that the definition of pre-opers and the formation of the filtration
\(\{V_i\}_{0\leq i\leq p}\) can be made on any smooth \(1\)-dimensional
noetherian scheme over \(k\), and statements (i) and (ii) make sense in this
context.  In fact, the statements being local, we are reduced to the
case of a free \(\ringO_X\)-module \(W\).  Moreover, all the relevant
formations commute with direct sums, hence we are reduced to the case
of \(W=\ringO_X\).

We can even reduce to the case \(X=\Spec k[[t]]\) and use an explicit
calculation to complete the proof.  Alternatively, one can check that
the construction of \cite[Remarkques 4.1.2 (2)]{raynaud82a}
gives the same filtration and proves the theorem.\qed

\begin{rem*}
  If \(L\) is a line bundle such that \(V=F_*L\) is stable, then the
  above result shows that \(F^*V\) is highly unstable, and it is
  likely that \(V\) is in a minimal Frobenius stratum (this is
  indeed the case in characteristic two).  At least, the bound in
  \cite[Theorem~3.1, page 51]{Su} is reached by \(V\) when the rank is
  \(p\).

  We refer to Theorem~\ref{gen-cons-thm} (ii) for conditions ensuring
  stability of \(F_*L\).  The following is a partial converse to the
  above theorem.
\end{rem*}

\begin{prop*}
Assume that \(F_*L\) is stable for any line bundle \(L\).
Let $\bigr(V,\nabla,\{V_i\}_{0\leq i\leq p})$
be a $p$-flat pre-oper with $\rank(V) = p$ and $V^{\nabla=0}$
stable. Then $V^{\nabla=0}\isom F_*(L)$ for a suitable line bundle
$L$ on $X$.
\end{prop*}

\begin{proof}
As $0\subset V_1\subset\cdots \subset V_{p-1}\subset V_p=V$ is a pre-oper of rank
$p$, $V/V_{p-1}=L$ is a line bundle.  The morphism
$F^*(V^{\nabla=0})=V\to V/V_{p-1}=L$ gives by adjunction a
non-zero morphism $V^{\nabla=0}\to F_*(L)$.  Since
these bundles are stable and of the same degree,
the map $V^{\nabla=0}\to F_*(L)$
is an isomorphism.
\end{proof}

\subsection{} The underlying bundle of a pre-oper
is typically unstable.  In some circumstances, the Frobenius
descent of a \(p\)-flat pre-oper is (semi)-stable.  This provides
a way of constructing Frobenius destabilized bundles in terms
of pre-opers.

\begin{prop*}
Let $(V,\nabla,\{V_i\}_{0\leq i\leq l})$ be a $p$-flat pre-oper with
\(V_1\) semi-stable of rank \(r_1\).
\begin{itemize}
\item [\rm (i)]
Suppose
\(
p > l^2(l-1)(g-1)r_1.
\)
Then $V^{\nabla=0}$ is semi-stable.
\item [\rm (ii)] Suppose \(r_1=1\) and \(p > l^2(l-1)(g-1)\).  Then
\(V^{\nabla=0}\) is stable.
\end{itemize}
\end{prop*}

\begin{proof}  (i) A direct computation shows that
\(\mu(V_1)=\mu(V)+(l-1)(g-1)/r_1\).  Since
  \(V_1\) is semi-stable, \(\{V_i\}_{0\leq i \leq l}\) is the
  Harder-Narasimhan filtration of \(V\).  In particular, if
  \(W\ \subset V\), then \(\mu(W) \leq \mu(V_1)\).

Let \(V'=V^{\nabla=0}\).  Then \(F^*(V')=V\).  Suppose that \(V'\) is not
semi-stable, and \(W' \subset V'\) is such that \(\mu(W')>\mu(V')\).
Then \(W=F^*(W')\) satisfies \(\mu(W)\leq \mu(V_1)\),
and
\[
\frac{\mu(V)}{p}=\mu(V')<\mu(W')=\frac{\mu(W)}{p}\leq
\frac{\mu(V_1)}{p}=\frac{\mu(V)}{p}+\frac{(l-1)(g-1)}{p \cdot r_1}.
\]
However, \(\mu(W')\) and \(\mu(V')\) are fractions of the form \(a/b\),
with \(a, b\in\bbZ\), \(0<b\leq l \cdot r_1\).  Therefore,
\(\mu(W')-\mu(V')\geq(l \cdot r_1)^{-2}\).
This contradicts the assumption on $p$.

(ii) Let \(V'=V^{\nabla=0}\) and \(0=W_0' \subset \cdots \subset
W_s'=V'\) be a Jordon-H\"older series for \(V'\).  Then each
\(W_i'/W'_{i-1}\) is stable of slope \(\mu/p\), where \(\mu=\mu(V)\).  Let
\(W_i=F^*(W_i'/W'_{i-1})\), \(\mu_{\rm max}(W_i)\) be the largest possible
slope of sub-bundles of \(W_i\) and \(\mu_{\rm min}(W_i)\) be the smallest
possible slope of quotient-bundles of \(W_i\).  By definition,
\(\mu_{\rm min}(W_i) \leq \mu \leq \mu_{\rm max}(W_i)\).

Let \(i_0\) be the smallest integer such that \(F^*(W'_{i_0})\to
V_l/V_{l-1}\) is
non-zero.  Then \(\mu_{\rm min}(W_{i_0})\leq\mu(V_l/V_{l-1})=\mu-(l-1)(g-1)\).  Similarly,
there exists an index \(i_1\) such that \(\mu_{\rm max}(W_{i_1}) \geq
\mu(V_1)=\mu+(l-1)(g-1)\).

A theorem of Sun \cite[Theorem 3.1]{Su} asserts that
\[
\mu_{\rm max}(W_i)-\mu_{\rm min}(W_i) \leq (\rank(W_i)-1)(2g-2).
\]
This implies that \(\rank(W_{i_0}) \ge (l+1)/2\) and
\(\rank(W_{i_1}) \ge (l+1)/2\).  Thus \(i_0=i_1\) and
\[
(l-1)(2g-2)\leq \mu_{\rm max}(W_{i_0})-\mu_{\rm min}(W_{i_0})\leq (\rank(W_{i_0})-1)(2g-2)
\]
by Sun's theorem again.
Therefore, \(\rank W_{i_0} \geq l\) and this forces
\(W'_{i_0}\) to be the only Jordan-H\"older factor of \(V'\).
\end{proof}

\begin{rem*} The bound on \(p\) can often be improved for particular
\((l,g,r_1,\deg(V_1))\).  This is clear from the proof.
\end{rem*}

\begin{center}\small
\begin{tabular}{l@{\hspace{0.6in}}l}
Kirti Joshi&S. Ramanan\\
Department of Mathematics&School of Mathematics\\
University of Arizona&Tata Institute of Fundamental Research\\
Tucson&Homi Bhabha Road\\
AZ 85721&Mumbai 400 005\\
USA                & India\\
\noalign{\smallbreak} {\it Email:\/}
{\tt kirti@math.arizona.edu}&
{\it Email:\/} {\tt ramanan@math.tifr.res.in}\\
\end{tabular}
\end{center}

\begin{center}\small
\begin{tabular}{l@{\hspace{0.6in}}l}
Eugene Z. Xia&Jiu-Kang Yu\\
National Center for Theoretical Sciences & National Center for Theoretical Sciences\\
National Tsing Hua University & National Tsing Hua University\\
No. 101, Sec 2, Kuang Fu Road, & No. 101, Sec 2, Kuang Fu Road,\\
Hsinchu, Taiwan 30043, Taiwan R.O.C. & Hsinchu, Taiwan 30043, Taiwan R.O.C.\\
\noalign{\smallbreak}
{\it Email:\/} {\tt xia@math.umass.edu}&
{\it Email:\/} {\tt yu@math.umd.edu}\\
\end{tabular}
\end{center}

\end{document}